\documentclass[11pt]{article}
\usepackage{amssymb}
\usepackage{amscd}
\usepackage{amsmath}
\usepackage{amsthm}

\let\geq\geqslant
\let\leq\leqslant

\renewcommand{\a}{{\alpha}}
\renewcommand{\b}{{\beta}}
\newcommand{\g}{{\gamma}}
\renewcommand{\d}{{\delta}}
\newcommand{\e}{{\varepsilon}}

\renewcommand{\l}{{\lambda}}

\def\o{{\omega}}

\newcommand{\s}{{\sigma}}

\newcommand{\z}{{\zeta}}
\newcommand{\D}{{\Delta}}
\newcommand{\G}{{\Gamma}}
\def\L{{\Lambda}}

\DeclareMathOperator{\sign}{sign}

\DeclareMathOperator{\supp}{supp}
\DeclareMathOperator{\rank}{rank}

\DeclareMathOperator{\Tr}{Tr}

\DeclareMathOperator{\diag}{diag}

\renewcommand\Im{\hbox{{\rm Im}}\,}
\renewcommand\Re{\hbox{{\rm Re}}\,}
\newcommand{\abs}[1]{\lvert#1\rvert}
\newcommand{\aabs}[1]{\left\lvert#1\right\rvert}
\newcommand{\norm}[1]{\lVert#1\rVert}
\newcommand{\nnorm}[1]{\left\lVert#1\right\rVert}
\newcommand{\br}[1]{\left(#1\right)}

\newcommand{\loc}{{\hbox{\rm\scriptsize loc}}}

\newcommand{\R}{{\mathbb R}}
\newcommand{\N}{{\mathbb N}}
\newcommand{\Z}{{\mathbb Z}}
\newcommand{\C}{{\mathbb C}}
\newcommand{\T}{{\mathbb T}}

\numberwithin{equation}{section}

\theoremstyle{plain}
\newtheorem{theorem}{\bf Theorem}[section]
\newtheorem{lemma}[theorem]{\bf Lemma}
\newtheorem{proposition}[theorem]{\bf Proposition}

\theoremstyle{definition}

\theoremstyle{remark}
\newtheorem*{remark*}{\bf Remark}
\newtheorem{remark}[theorem]{\bf Remark}

\renewcommand{\qed}{\vrule height7pt width5pt depth0pt}

\newcommand{\B}{{B}}
\renewcommand{\AA}{{\mathcal A}}
\newcommand{\BB}{{\mathcal B}}
\newcommand{\M}{{\mathcal M}}

\newcommand{\wt}{\widetilde}

\newcommand{\J}{{\mathbb J}}
\renewcommand{\aa}{\kappa}
\newcommand{\oo}{{\vec{\omega}}}
\newcommand{\vp}{{\rm v.p.\,}}

\newcommand{\mom}{{\mu}}

\begin{document}

\title{A trace formula and high energy spectral asymptotics for the perturbed Landau Hamiltonian}

\date{2 October 2003}

\author{E.~Korotyaev\thanks{
Institut f\"ur Mathematik, Humboldt Universit\"at zu Berlin,
Rudower Chaussee 25, 12489, Berlin, Germany.
e-mail: ek@mathematik.hu-berlin.de}\mbox{ } and
A.~Pushnitski\thanks{
Department of Mathematical Sciences, Loughborough University,
Loughborough, LE11 3TU, U.K.
e-mail: a.b.pushnitski@lboro.ac.uk}}

\maketitle

\begin{abstract}
A two-dimensional Schr\"odinger operator with a constant
magnetic field perturbed by a smooth compactly supported potential
is considered. The spectrum of this operator consists of
eigenvalues which accumulate to the Landau levels.
We call the set of eigenvalues near the $n$'th Landau level
an $n$'th eigenvalue cluster, and study the distribution of
eigenvalues in the $n$'th cluster as $n\to\infty$.
A complete asymptotic expansion for the eigenvalue moments
in the $n$'th cluster is obtained and some coefficients
of this expansion are computed. A  trace formula involving the
first eigenvalue moments is obtained.
\end{abstract}

\noindent

\bigskip

\section{Introduction and Main Results}
\label{sec.a}

\textbf{1. Introduction}
Let $H$ in $L^2(\R^2,dx_1\,dx_2)$ be the following magnetic Schr\"odinger
operator:
$$
H=\biggl(-i\frac{\partial}{\partial x_1}+\frac{\B}{2}x_2\biggr)^2 +
\biggl(-i\frac{\partial}{\partial x_2}-\frac{\B}{2}x_1\biggr)^2, \quad \B>0.
$$
The operator $H$ describes a quantum particle in $\R^2$ in a constant
homogeneous magnetic field of the magnitude $\B$;
it is often called the Landau Hamiltonian.
It is well known
\cite{LL}
that the spectrum of $H$ consists of a sequence of eigenvalues (Landau levels)
$\L_n=\B(2n+1)$, $n\in\Z_+\equiv\{0,1,2,\dots\}$. Each of these eigenvalues has
infinite multiplicity.

Let $V\in C_0^\infty(\R^2)$ be a real valued function (potential in physical terminology).
Consider the spectrum of the operator $H+V$.
It is well known (see \cite{AHS}) that $V$ is a relatively
compact perturbation of $H$  and therefore the essential
spectrum of $H+V$ is the same as that of $H$, i.e. consists of the Landau levels.
The operator $H+V$ may have eigenvalues of finite multiplicities
which can accumulate to the Landau levels.

Let us define disjoint intervals $\D_0=[\inf\s(H+V),2B)$,
$\D_n=[\L_n-B,\L_n+B)$,
$n\in\N$, so that $\s(H+V)\subset \cup_{n=0}^\infty\D_n$.
We shall call the set $\s(H+V)\cap\D_n$ the $n$'th \emph{eigenvalue cluster}.
For a fixed $n$, the distribution of eigenvalues in the $n$'th cluster was studied in
\cite{MR}, \cite{RW} (these papers contain also references to earlier work
on this problem).
It was found that eigenvalues accummulate to $\L_n$ exponentially fast.

Our aim is to study the asymptotic distribution of eigenvalues in the $n$'th cluster
as $n\to\infty$.
Our first preliminary result is that the width of the $n$'th cluster is $O(n^{-1/2})$:

\begin{proposition}\label{thm.a.1}
There exist $C>0$ and $N\in\N$ such that for all $n\geq N$, one has
$$
\s(H+V)\cap\D_n\subset(\L_n-Cn^{-1/2},\L_n+Cn^{-1/2}).
$$
The constants $C$ and $N$ depend only on
$\sup\limits_{x\in\R^2} \abs{V(x)}$ and on the diameter of $\supp V$.
\end{proposition}
The power $n^{-1/2}$ in the above Proposition is sharp;
see Remark~\ref{rmk.d.2} below.

\textbf{2. Definition of eigenvalue moments $\mom_n$.}
We would like to define moments of eigenvalues in the $n$'th cluster.
First, in order to explain the main idea of the definition,
let us define the eigenvalue moments `naively'
for the case $\norm{V}<\B$;
here and in what follows $\norm{V}\equiv\norm{V}_{L^\infty}$.
Fix $n\in\Z_+$ and enumerate $\l_1$, $\l_2$, $\l_3$,\dots all eigenvalues
in the $n$'th cluster so that
$\abs{\l_1-\L_n}\geq \abs{\l_2-\L_n}\geq \abs{\l_3-\L_n}\geq\cdots$.
Let us define the eigenvalue moments by
\begin{equation}
\label{a.1}
\mom_n=\sum_j (\l_j-\L_n),
\quad n\in\Z_+
\quad (\norm{V}<\B).
\end{equation}
By the above quoted result of \cite{MR}, \cite{RW}, the
rate of convergence $\l_j\to\L_n$ as $j\to\infty$ is
exponential, and therefore the series \eqref{a.1}
converges absolutely.

In order to give the definition of eigenvalue moments which is suitable
both for the case $\norm{V}<\B$ and for the case $\norm{V}\geq \B$,
we need to recall the notion of the spectral shift function
for the pair of operators $H+V$, $H$. The spectral shift function was introduced in
an abstract operator theoretic setting in
\cite{L,K}; see also the book \cite{Ya}.
Recall that (see \cite{AHS})
\begin{equation}
(H+V-\l_0)^{-1}-(H-\l_0)^{-1}\in\text{Trace class}, \quad \l_0<\inf\s(H+V).
\label{a.2}
\end{equation}
This enables one to define the spectral shift function $\xi\in L^1_\loc(\R)$
for the pair $H+V$, $H$.
The spectral shift function $\xi$ is determined by the following two conditions:

(i) For any `test function' $\phi\in C_0^\infty(\R)$,
one has the trace formula:
\begin{equation}
\Tr(\phi(H+V)-\phi(H))=
\int_{-\infty}^\infty \xi(\l)\phi'(\l)d\l. \label{a.3}
\end{equation}

(ii) $\xi(\l)=0$ for $\l<\inf\s(H+V)$.\newline
Note that $\phi(H+V)-\phi(H)\in\text{Trace class}$ by \eqref{a.2} (see \cite{Ya}).
In fact, the class of admissible test functions $\phi$ is much wider than
$C_0^\infty(\R)$. In particular, this class includes exponentials
$\phi(\l)=e^{-t\l}$, $t>0$; we will use the latter fact in the sequel.

Condition (i) determines the spectral shift function up to an additive constant;
condition (ii) fixes this constant.
As it follows from the trace formula \eqref{a.3},
for $\l\in\R\setminus\s(H+V)$
the spectral shift function can be determined by
(see \cite{Ya}, section 8.7 and formula (8.2.20))
\begin{equation}
\xi(\l)=\Tr(E_H(\l)- E_{H+V}(\l)),
\quad \l\in\R\setminus\s(H+V),
\label{a.4}
\end{equation}
where $E_H(\l)$ and $E_{H+V}(\l)$ are the spectral projections of $H$ and $H+V$
associated with the interval $(-\infty,\l)$.
In particular, it follows that
$\xi$ is constant and integer-valued on the intervals of the set $\R\setminus\s(H+V)$.

Now we are ready to give a general definition of the eigenvalue
moments:
\begin{equation}
\mom_n=
\int_{\D_n}\xi(\l)d\l,
\quad n\in\Z_+.
\label{a.5a}
\end{equation}
Let us explain why
the definitions \eqref{a.5a} and \eqref{a.1} coincide for $\norm{V}<\B$.
{}From \eqref{a.4} one can see that for $\norm{V}<\B$
\begin{equation}
\xi(\l)=
\begin{cases}
&\text{the number of eigenvalues of $H+V$ in $(\l,\L_n+\B)$
if $\l\in(\L_n,\L_n+\B)$;}
\\
&(-1)\times\text{the number of eigenvalues of $H+V$ in $(\L_n-\B,\l)$
if $\l\in(\L_n-\B,\L_n)$.}
\end{cases}
\label{a.6}
\end{equation}
{}From here it follows that \eqref{a.5a} and \eqref{a.1} coincide.
\begin{remark*}
Proposition~\ref{thm.a.1} shows that
$$
\s(H+\tau V)\cap\D_n\subset(\L_n-Cn^{-1/2},\L_n+Cn^{-1/2}),
\quad
\forall \tau \in[0,1]
$$
for any $V$ and
all sufficiently large $n$. From here, using \eqref{a.4} and a continuity
in $t$ argument, one can prove that for any $V$ and all sufficiently large $n$,
\begin{equation}
\supp \xi\cap\D_n\subset(\L_n-Cn^{-1/2},\L_n+Cn^{-1/2})
\label{d.11}
\end{equation}
and \eqref{a.6} holds true.
Thus,  definition \eqref{a.1} is applicable for any $V$ and all
sufficiently large $n$.
\end{remark*}

\textbf{3. Main result}

\begin{theorem}\label{thm.a.2}
The asymptotic expansion
\begin{equation}
\mom_n\sim
\a_0+\frac{\a_1}{n^{1/2}}+\frac{\a_2}{n}+ \frac{\a_3}{n^{3/2}}+\cdots,
\quad n\to\infty,
\label{a.7a}
\end{equation}
holds true with some real coefficients $\a_j$. Moreover, one has
\begin{equation}
\a_0=\frac{\B}{2\pi}\int_{\R^2} V(x)dx,\quad
\a_1=\a_2=0,\quad
\a_3=-\frac{\sqrt{\B}}{16\sqrt{2}\pi^3}
\int_{\R^2}\int_{\R^2} \frac{V(x)V(y)}{\abs{x-y}}dxdy.
\label{a.8}
\end{equation}
The identity (trace formula)
\begin{equation}
\sum_{n=0}^\infty\biggl(
\mom_n-\frac{\B}{2\pi}\int_{\R^2} V(x)dx\biggr)
=-\frac1{8\pi}\int_{\R^2} V^2(x)dx
\label{a.10}
\end{equation}
holds true.
\end{theorem}
\textbf{Remarks}
(1)
The coefficients $\a_0$, $\a_1$, $\a_2$ are obtained by comparing
the asymptotic expansion \eqref{a.7a} with the small $t$ asymptotic
expansion of $\Tr(e^{-t(H+V)}-e^{-tH})$ --- see section~\ref{sec.b} below.
It does not seem possible to obtain the coefficient $\a_3$ by using a similar
argument; we obtain it by a more direct analysis (see end of section~\ref{sec.e}).
(2)
Some results concerning the eigenvalue distribution in clusters
for large $n$ can be found in \cite{Ta}.
(3)
Similar trace formula for the two-dimensional perturbed harmonic
oscillator was obtained in \cite{FM}.
(4)
It might be interesting to note in connection with the formula for
$\a_3$ that the integral $\int\int\frac{V(x)V(y)}{\abs{x-y}}dxdy$
appears as the coefficient of the leading term in the high energy
asymptotic expansion of the total scattering cross-section
for the pair of operators $-\Delta$, $-\Delta+V(x)$ in $L^2(\R^2)$.

Along with the moments $\mom_n$, we will use the
higher order moments
\begin{equation}
\mom_n^{(k)}=
(k+1)\int_{\D_n}(\l-\L_n)^k\xi(\l)d\l,
\quad k\in\N,\quad n\in\Z_+.
\label{a.5}
\end{equation}
In the case $\norm{V}<\B$, the last definition becomes
$\mom_n^{(k)}=\sum_j(\l_j-\L_n)^{k+1}$.
We will also prove the asymptotic expansion
\begin{equation}
\mom_n^{(k)}\sim n^{-k/2}
\bigl(\a^{(k)}_0+\frac{\a^{(k)}_1}{n^{1/2}}+
\frac{\a^{(k)}_2}{n}+ \frac{\a^{(k)}_3}{n^{3/2}}+\cdots\bigr),
\quad k\in\N,
\quad n\to\infty.
\label{a.7}
\end{equation}
Below for consistency we write $\mom_n\equiv\mom_n^{(0)}$,
$\a_j\equiv\a_j^{(0)}$.

\textbf{4.  The structure of the paper}
We will use three distinct arguments to
prove Proposition~\ref{thm.a.1}, the asymptotic expansions
\eqref{a.7a}, \eqref{a.7} and the trace formula \eqref{a.10}.
The proof of Proposition~\ref{thm.a.1} is self-contained,
the proof of the asymptotic expansions \eqref{a.7a}, \eqref{a.7}
depends on the estimates \eqref{d.2}, \eqref{d.3}
obtained in the proof of
Proposition~\ref{thm.a.1}, and the proof
of the trace formula \eqref{a.10} depends on both
Proposition~\ref{thm.a.1} and the expansions \eqref{a.7a}, \eqref{a.7}.

In section~\ref{sec.b}, assuming Proposition~\ref{thm.a.1}
and the existence of the asymptotic expansions \eqref{a.7a}, \eqref{a.7}, we
prove the trace formula \eqref{a.10} and
derive formulae \eqref{a.8} for the coefficients $\a_0$, $\a_1$
and $\a_2$.
The argument is quite elementary.

Proposition~\ref{thm.a.1} is proven in section~\ref{sec.d}.

In Sections \ref{sec.e}, \ref{sec.f}, \ref{sec.fa} we justify the asymptotic
expansions \eqref{a.7a}, \eqref{a.7}. The proof is based on a detailed analysis of the
properties of the integral kernel of the resolvent of $H$
(see \eqref{f.1}) and on some  facts from the theory of confluent
hypergeometric functions.
This part of the paper is fairly elementary in nature but technically
is rather complicated.

In Section~\ref{sec.e} we also prove the formula \eqref{a.8}
for the coefficient  $\a_3$.

\textbf{5. Notation}
We use notation $\norm{A}$, $\norm{A}_{S_2}$, $\norm{A}_{S_1}$
for the operator norm, the Hilbert-Schmidt norm, and the trace class
norm of an operator $A$.
By $C$, $c$ we denote various constants in the estimates.

\section{Proof of the trace formula}
\label{sec.b}

\textbf{1. Heat kernel asymptotics.}
\begin{lemma}\label{l.3.1}
The asymptotic formula
\begin{equation}
\Tr(e^{-tH}- e^{-t(H+V)})=
\frac1{4\pi}\int_{\R^2} V(x)dx
-\frac{t}{8\pi} \int_{\R^2} V^2(x)dx+O(t^2),\quad  t\to+0
\label{b.20}
\end{equation}
holds true.
\end{lemma}
\begin{proof}
The required asymptotics can be obtained by using general results
on asymptotic expansions of heat kernels of second
order elliptic operators.
However, the two term asymptotic formula \eqref{b.20}
is considerably simpler than the aforementioned general
results, and so we prefer to give a direct `elementary' proof.
We use the formula
\begin{equation}
e^{-tH}-e^{-t(H+V)}=
\int_0^t e^{-(t-t_1)H}Ve^{-t_1 (H+V)}dt_1
\label{b.21}
\end{equation}
and the explicit formula for the integral kernel of $e^{-tH}$
(see \cite{AHS}):
\begin{equation}
e^{-tH}(x,y)=\frac{\B}{4\pi\sinh(\B t)}
\exp(-\tfrac{\B}{4}\abs{x-y}^2\coth(\B t)+i\tfrac{\B}{2}[x,y]),
\qquad x,y\in\R^2,\quad t>0,
\label{f.2a}
\end{equation}
where $[x,y]\equiv x_1y_2-x_2y_1$.
Iterating \eqref{b.21}, we obtain
\begin{align*}
\Tr(e^{-tH}&- e^{-t(H+V)})=I_1(t)+I_2(t)+I_3(t),
\\
I_1(t)&=\int_0^t\Tr (e^{-(t-t_1)H}Ve^{-t_1 H})dt_1=
t\Tr(Ve^{-tH})=\frac{1}{4\pi}\int_{\R^2}V(x)dx+O(t^2),
\quad t\to+0,
\\
I_2(t)&=\int_0^tdt_1\int_0^{t_1}dt_2
\Tr(e^{-(t-t_1)H}Ve^{-(t_1-t_2)H}Ve^{-t_2H}),
\\
I_3(t)&=\int_0^t dt_1\int_0^{t_1}dt_2 \int_0^{t_2} dt_3
\Tr(e^{-(t-t_1)H}Ve^{-(t_1-t_2)H}Ve^{-(t_2-t_3)H}Ve^{-t_3(H+V)}).
\end{align*}
Consider the term $I_2(t)$. Introducing the new variable
$s=t-t_1+t_2$, we have
\begin{multline*}
I_2(t)=\int_0^tds\, s \Tr (e^{-sH}Ve^{-(t-s)H}V)
\\
=\int_0^tds\,\frac{s}{\sinh \B s\sinh \B(t-s)}
\int_{\R^2}dz\,
\exp\br{-\tfrac{\B}{4}(\coth \B s+\coth\B(t-s))\abs{z}^2}
W(z),
\end{multline*}
where $W(z)=\frac{B^2}{(4\pi)^2}\int_{\R^2} V(y)V(y+z)dy$.
Let us write $W(z)=W(0)+(W(z)-W(0))$ and split $I_2(t)$
accordingly: $I_2(t)=I_2^{(1)}(t)+I_2^{(2)}(t)$.
For $I_2^{(1)}(t)$, explicitly computing the integrals,
we obtain
$$
I_2^{(1)}(t)=W(0)\frac{2\pi t^2}{\B\sinh\B t}
=
\frac{t}{8\pi}\int_{\R^2}V^2(x)dx+O(t^3),
\quad t\to+0.
$$
For $I_2^{(2)}(t)$, using the estimate $\abs{W(z)-W(0)}\leq
C\abs{z}$, we get
\begin{multline*}
\abs{I_2^{(2)}(t)}\leq
C\int_0^tds\,\frac{s}{\sinh\B s\sinh\B(t-s)}
\int dz\,\abs{z}\exp(-\tfrac{\B}{4}(\coth \B s+\coth\B(t-s))\abs{z}^2)
\\
=
\frac{C_1}{\sinh\B t}\int_0^t ds\,
s(\sinh \B s\sinh \B(t-s))^{1/2}
=O(t^2)\quad t\to+0.
\end{multline*}
Finally, let us estimate the term $I_3(t)$.
Using the Hilbert-Schmidt norm estimate
$\norm{e^{-tH}V}_{S_2}\leq C t^{-1/2}$, $t>0$,
we obtain
\begin{multline*}
\abs{\Tr(e^{-(t-t_1)H}Ve^{-(t_1-t_2)H}Ve^{-(t_2-t_3)H}Ve^{-t_3(H+V)})}
\\
\leq
\norm{e^{-(t-t_1)H}V}_{S_2}\norm{e^{-(t_1-t_2)H}V}_{S_2}
\norm{V}\norm{e^{-t(H+V)}}
\leq \frac{C}{\sqrt{t-t_1}\sqrt{t_1-t_2}},
\end{multline*}
which yields $I_3(t)=O(t^2)$, $t\to+0$.
\end{proof}

\textbf{2. Auxiliary estimate}
\begin{lemma}\label{lma.b.2}
One has
$$
\int_{\Delta_n}\abs{\xi(\l)}d\l=O(1),
\quad n\to\infty.
$$
\end{lemma}

\begin{proof}
Recall that the spectral shift function is monotone with respect to the perturbation $V$.
I.e., denoting temporarily by $\xi(\l;V)$ the spectral shift function
corresponding to the potential $V$, we have
$$
\text{ if  $V_1\leq V_2$, then }
\xi(\l;V_1)\leq \xi(\l;V_2)
\text{ a.e. }\l\in\R.
$$
Also, we have $\xi(\l;-V)=-\xi(\l;V)$.
Let us choose $V_1,V_2\in C_0^\infty(\R^2)$ such that
$$
V_1\geq0,\quad V_2\geq0\quad\text{and} \quad -V_2\leq V\leq V_1.
$$
Then
$$
-\xi(\l;V_2)=\xi(\l;-V_2)\leq\xi(\l;V)\leq\xi(\l;V_1),
\quad \xi(\l;V_1)\geq0, \quad \xi(\l;V_2)\geq0.
$$
Therefore, $\abs{\xi(\l;V)}\leq \xi(\l;V_1)+\xi(\l;V_2)$ and
$$
\int_{\D_n}\abs{\xi(\l;V)}d\l\leq
\int_{\D_n}\xi(\l;V_1)d\l+\int_{\D_n}\xi(\l;V_2)d\l.
$$
By the asymptotics \eqref{a.7} for $\mom^{(0)}_n$, applied to the
potentials $V_1$ and $V_2$, the r.h.s. of the last inequality is
$O(1)$ as $n\to\infty$.
\end{proof}

\textbf{3. Proof of the trace formula \eqref{a.10}
and formulae \eqref{a.8} for $\a^{(0)}_0$, $\a^{(0)}_1$, $\a^{(0)}_2$.    }

1. By  Krein's trace formula \eqref{a.3} with $\phi(\l)=e^{-t\l}$, $t>0$,
and Lemma~\ref{l.3.1}, we obtain
\begin{equation}
\int_{-\infty}^\infty \xi(\l) e^{-t\l}d\l= \frac1t \Tr(e^{-tH}-
e^{-t(H+V)})= \frac1{4\pi t}\int_{\R^2} V(x)dx -\frac1{8\pi}
\int_{\R^2} V^2(x)dx+O(t), \label{b.3}
\end{equation}
as $t\to+0$. On the other hand, one can rewrite the integral in
the l.h.s. of \eqref{b.3} as a sum over the eigenvalue clusters:
\begin{equation}
\int_{-\infty}^\infty \xi(\l) e^{-t\l}d\l= \sum_{n=0}^\infty
\int_{\D_n} \xi(\l) e^{-t\l}d\l.
\label{b.3a}
\end{equation}
Let us use Taylor's  formula for $e^{-t\l}$, $\l\in\D_n$:
\begin{equation}
\abs{e^{-t\l}-e^{-t\L_n}(1-t(\l-\L_n)+\frac12 t^2(\l-\L_n)^2)}
\leq
C t^3\abs{\l-\L_n}^3,\quad \l\in\D_n.
\label{b.3b}
\end{equation}
By \eqref{d.11} and Lemma \ref{lma.b.2}, we have
\begin{equation}
\sum_{n=0}^\infty \int_{\D_n}\abs{\xi(\l)} \abs{\l-\L_n}^3d\l=
\sum_{n=0}^\infty O(n^{-3/2})<\infty.
\label{b.3c}
\end{equation}
Combining \eqref{b.3a}--\eqref{b.3c}, we obtain
\begin{equation}
\int_{-\infty}^\infty \xi(\l) e^{-t\l}d\l= \sum_{n=0}^\infty \mom^{(0)}_n e^{-t\L_n}
-t\sum_{n=0}^\infty \mom^{(1)}_n e^{-t\L_n} +\frac12 t^2\sum_{n=0}^\infty
\mom^{(2)}_n e^{-t\L_n} +O(t^3),\quad
t\to+0.
\label{b.4}
\end{equation}
Below we compare \eqref{b.3} and \eqref{b.4}.

2. We will use the following elementary formulae for $t\to+0$:
\begin{align}
&\sum_{n=0}^\infty e^{-t\L_n}=
\frac{e^{-Bt}}{1-e^{-2Bt}}=
\frac1{2Bt}+O(t);
\label{b.5}
\\
&\sum_{n=1}^\infty n^{-1/2} e^{-t\L_n}=
\int_0^\infty x^{-1/2} e^{-tB(2x+1)}dx+O(1)= \frac{\sqrt{\pi}}{\sqrt{2Bt}}+O(1);
\label{b.6}
\\
&\sum_{n=1}^\infty n^{-1} e^{-t\L_n}=
\int_1^\infty x^{-1} e^{-tB(2x+1)}dx+O(1)=
-\log t +O(1);
\label{b.7}
\\
&\sum_{n=1}^\infty n^{-3/2} (1-e^{-t\L_n})=
\int_0^t \bigl(\sum_{n=1}^\infty n^{-3/2}\L_n e^{-s\L_n}\bigr)ds
=2\sqrt{2\pi \B t}+O(t).
\label{b.7a}
\end{align}
Using \eqref{b.5}--\eqref{b.7} and the asymptotic expansions
\eqref{a.7a}, \eqref{a.7} for $\mom_n^{(k)}$, we obtain for $t\to+0$:
\begin{align}
&\sum_{n=0}^\infty \mom_n^{(0)} e^{-t\L_n}= \frac{\a_0^{(0)}}{2Bt}+
\frac{\a_1^{(0)}\sqrt{\pi}}{\sqrt{2Bt}}- \a_2^{(0)}\log t+O(1);
\label{b.8}
\\
&\sum_{n=0}^\infty \mom_n^{(1)} e^{-t\L_n}=
\a_0^{(1)}\frac{\sqrt{\pi}}{\sqrt{2Bt}}+ O(\log t);
\label{b.9}
\\
&\sum_{n=0}^\infty \mom_n^{(2)} e^{-t\L_n}= O(\log t).
\notag
\end{align}
Substituting this into \eqref{b.4} and comparing with \eqref{b.3}, we find:
$$
\a^{(0)}_0=\frac{B}{2\pi}\int_{\R^2} V(x)dx,\quad \a^{(0)}_1=\a^{(0)}_2=0.
$$
Thus, $\mom_n^{(0)}=\a_0^{(0)}+\a_3^{(0)}n^{-3/2}
+O(n^{-2})$, and so we get:
\begin{multline}
\sum_{n=0}^\infty \mom_n^{(0)} e^{-t\L_n}
=
\a_0^{(0)}\sum_{n=0}^\infty e^{-t\L_n}
+\sum_{n=0}^\infty (\mom_n^{(0)}-\a_0^{(0)})
+\sum_{n=1}^\infty (\mom_n^{(0)}-\a_0^{(0)})(e^{-t\L_n}-1)
\\
=
\frac{\a_0^{(0)}}{2Bt}+
\sum_{n=0}^\infty (\mom_n^{(0)}-\a_0^{(0)})
+\a_3^{(0)}\sum_{n=1}^\infty n^{-3/2}(e^{-t\L_n}-1)
+\sum_{n=1}^\infty (\mom_n^{(0)}-\a_0^{(0)}-\a_3^{(0)}n^{-3/2})(e^{-t\L_n}-1)
+O(t)
\\
=
\frac{\a_0^{(0)}}{2Bt}+
\sum_{n=0}^\infty (\mom_n^{(0)}-\a_0^{(0)})
-\a_3^{(0)}2\sqrt{2\pi Bt}+o(\sqrt{t}),
\label{b.10}
\end{multline}
as $t\to+0$.
Upon comparing with \eqref{b.3}, we find the trace formula \eqref{a.10}
and also the formula
\begin{equation}
\a^{(1)}_0=-4\B\a^{(0)}_3.
\label{a.9}
\end{equation}
We will use \eqref{a.9} later in section~\ref{sec.e}
in order to determine the coefficient $\a^{(0)}_3$.
\qed

\section{Proof of Proposition~\protect\ref{thm.a.1}}
\label{sec.d}
Let $P_n$, $n\geq0$, be the orthogonal projection onto the eigenspace of
the operator $H$ corresponding to the Landau level $\L_n$.
An explicit formula for the integral kernel of $P_n$ is available
(see e.g. \cite{RW}):
\begin{equation}
\label{d.1}
P_n(x,y)=\frac{\B}{2\pi}
L_n(\tfrac{\B}{2}\abs{x-y}^2)\exp(-\tfrac\B4\abs{x-y}^2+i\tfrac\B2[x,y]),
\quad x,y\in\R^2,
\end{equation}
where $L_n$ is the Laguerre polynomial and $[x,y]=x_1 y_2-x_2 y_1$.
\begin{lemma}\label{l.d.1}
Let $V$ be any bounded function on $\R^2$
with compact support. Then
\begin{align}
\label{d.2}
\norm{\abs{V}^{1/2}P_n\abs{V}^{1/2}}&=O(n^{-1/2}),
\quad n\to\infty,
\\
\label{d.3}
\norm{\abs{V}^{1/2}P_n\abs{V}^{1/2}}_{S_2}&=O(n^{-1/4}),
\quad n\to\infty.
\end{align}
\end{lemma}
In the proof of Proposition~\ref{thm.a.1},
we will only need the operator norm
estimate \eqref{d.2}.
The Hilbert-Schmidt norm estimate \eqref{d.3}
will be used in section \ref{sec.e}.
\begin{proof}
We will use the following asymptotic formula for the Laguerre polynomials
\cite[10.15(2)]{BE}:
\begin{equation}
\label{d.4}
L_n(t)=e^{t/2}J_0(\sqrt{(4n+2)t})+R_n(t),
\quad
R_n(t)=O(n^{-3/4}),\quad n\to\infty,
\end{equation}
where the bound $O(n^{-3/4})$
is uniform on any bounded sub-interval of $[0,\infty)$.
Let us write
$$
\abs{V}^{1/2}P_n\abs{V}^{1/2}=\AA_n+\BB_n,
$$
where $\AA_n$ and $\BB_n$ are the operators in $L^2(\R^2)$
with the integral kernels
\begin{align}
\label{d.6}
\AA_n(x,y)&=\frac{\B}{2\pi}J_0(\sqrt{\L_n}\abs{x-y})
e^{i\frac{\B}{2}[x,y]}\abs{V(x)}^{1/2}\abs{V(y)}^{1/2},
\\
\BB_n(x,y)&=\frac{\B}{2\pi}R_n(\tfrac{\B}{2}\abs{x-y}^2)
e^{i\frac{\B}{2}[x,y]}\abs{V(x)}^{1/2}\abs{V(y)}^{1/2}.
\notag
\end{align}
As $V$ is bounded and compactly supported,  \eqref{d.4}
gives
\begin{equation}
\label{d.7}
\norm{\BB_n}\leq \norm{\BB_n}_{S_2}
=O(n^{-3/4}),
\quad n\to\infty.
\end{equation}
Let us prove the estimate \eqref{d.3}.
By \eqref{d.7}, we only need to check that
$$
\norm{\AA_n}_{S_2}=O(n^{-1/4}).
$$
The latter estimate immediately follows from
the explicit form of the kernel of $\AA_n$
and from the simple inequality
$\abs{J_0(t)}\leq C/\sqrt{t}$ for $t>0$.

Next, let us prove the estimate \eqref{d.2}.
Let  $\wt \AA_n$ be the operator in $L^2(\R^2)$
with the integral kernel
$$
\wt \AA_n(x,y)=\frac{\B}{2\pi}J_0(\sqrt{\L_n}\abs{x-y})
\abs{V(x)}^{1/2}\abs{V(y)}^{1/2}.
$$
Note that up to a constant,
$\wt \AA_n$ coincides with the imaginary part of the sandwiched
resolvent of the operator $-\Delta$ in $L^2(\R^2)$:
$$
\wt \AA_n=\frac{2\B}{\pi}\Im(\abs{V}^{1/2}(-\Delta-\L_n-i0)^{-1}\abs{V}^{1/2}).
$$
It is well known (see \cite{A}) that
$$
\norm{\abs{V}^{1/2}(-\Delta-\l-i0)^{-1}\abs{V}^{1/2}}=O(\l^{-1/2}),
\quad \l\to\infty.
$$
Thus, we obtain
\begin{equation}
\norm{\wt \AA_n}=O(n^{-1/2}),\quad n\to\infty.
\label{d.8}
\end{equation}
Next, observe that the kernel of $\wt \AA_n$ differs from that of $\AA_n$
by a factor $e^{i\frac{\B}{2}[x,y]}$.
We are going to use this observation and apply the theory of
`multipliers of kernels of integral operators' \cite{BS}.
Let $\Omega$ be a sufficiently large ball in $\R^2$ so that
$\supp V\subset \Omega$ and let $\rho\in L^\infty(\Omega\times\Omega)$.
For a Hilbert-Schmidt class operator $T$ on $L^2(\Omega)$ with the integral kernel
$T(\cdot,\cdot)\in L^2(\Omega\times \Omega)$, let $\wt T$ be the
operator with the integral kernel $T(x,y)\rho(x,y)$.
Evidently, $\wt T$ is also a Hilbert-Schmidt class operator and
one has the estimate
$\norm{\wt T}_{S_2}\leq \norm{\rho}_{L^\infty}\norm{T}_{S_2}$.

Next, suppose that the mapping $T\mapsto \wt T$ sends the trace class $S_1$
into itself and there is a trace class norm bound
$\norm{\wt T}_{S_1}\leq C(\rho) \norm{T}_{S_1}$.
Then, by duality between the trace class $S_1$ and the class
${\mathbb B}(L^2(\Omega))$ of all bounded operators on
$L^2(\Omega)$, the mapping $T\mapsto \wt T$ can be extended onto
${\mathbb B}(L^2(\Omega))$ and the norm bound
$\norm{\wt T}\leq C(\rho) \norm{T}$ holds true.
In this case $\rho$ is called a bounded multiplier on the class
${\mathbb B}(L^2(\Omega))$.

A  sufficient condition (see \cite{BS})
for $\rho$ to be a bounded multiplier on ${\mathbb B}(L^2(\Omega))$ is
$$
\sup_{x\in\Omega}\norm{\rho(x,\cdot)}_{H^s(\Omega)}<\infty,
\quad s>1,
$$
where $H^s(\Omega)$ is the standard Sobolev class.
Clearly, $\rho(x,y)=e^{i\frac{\B}{2}[x,y]}$ satisfies the above condition,
and therefore
$$
\norm{\AA_n}\leq C\norm{\wt \AA_n}=O(n^{-1/2}),\quad
n\to\infty,
$$
which yields \eqref{d.2}.
\end{proof}
\begin{proof}[\textbf{Proof of Proposition~\ref{thm.a.1}}]
The proof is valid for any bounded compactly supported potential.
By the Birman-Schwinger principle, it suffices to show that
for some $C>0$ and all sufficiently large $n$,
\begin{equation}
\norm{\abs{V}^{1/2}R(\l)\abs{V}^{1/2}}<1,
\quad\text{ for all }
\l\in\D_n,
\quad
\abs{\l-\L_n}>\frac{C}{\sqrt{n}},
\label{d.10}
\end{equation}
where $R(\l)=(H-\l)^{-1}$.
Choose $l\in\N$ sufficiently large so that
$\norm{V}/\L_l<1/2$,
and write $R(\l)$ as
$$
R(\l)=\sum_{k=n-l}^{n+l}\frac{P_k}{\L_k-\l}+\wt R(\l).
$$
Then, for $\l\in\D_n$,
$$
\norm{\abs{V}^{1/2}R(\l)\abs{V}^{1/2}}
\leq
\sum_{k=n-l}^{n+l}
\frac{\norm{\abs{V}^{1/2}P_k\abs{V}^{1/2}}}{\abs{\L_k-\l}}
+
\norm{\abs{V}^{1/2}\wt R(\l)\abs{V}^{1/2}}.
$$
By the choice of $l$, one has $\norm{\abs{V}^{1/2}\wt R(\l)\abs{V}^{1/2}}<1/2$.
On the other hand, by Lemma~\ref{l.d.1},
$$
\sum_{k=n-l}^{n+l}
\frac{\norm{\abs{V}^{1/2}P_k\abs{V}^{1/2}}}{\abs{\L_k-\l}}
\leq(2l+1)O(n^{-1/2})\max_{n-l\leq k\leq n+l}\abs{\L_k-\l}^{-1}
=O(n^{-1/2})\abs{\L_n-\l}^{-1}.
$$
Thus, we get \eqref{d.10} for sufficiently large $C>0$.
\end{proof}
\begin{remark}\label{rmk.d.2}
1. The operator norm estimate \eqref{d.2} is sharp, i.e. for
any $V$ not identically zero, one can prove that
$$
\norm{\abs{V}^{1/2}P_{n}\abs{V}^{1/2}}\geq c/\sqrt{n}
$$
for some $c>0$ and all sufficiently large $n$.
Indeed, without the loss of generality assume that $\supp V$ contains an open
neighbourhood of zero and let $\wt V$ be a spherically symmetric potential,
$\wt V(x)=v(\abs{x})\leq \abs{V(x)}$.
Then by Lemma~3.3 of \cite{RW}, the eigenvalues $\l_{k}$ of
$\abs{\wt V}^{1/2}P_n\abs{\wt V}^{1/2}$ are given by
\begin{equation}
\frac{n!}{(n+k)!}\int_0^\infty v(\sqrt{2t/\B})e^{-t}t^k L_n^{(k)}(t)^2dt,
\quad k=-n,-n+1,-n+2,\dots,
\label{d.9}
\end{equation}
where $L_n^{(k)}$ are the Laguerre polynomials.
Taking $k=0$ and using the asymptotics \eqref{d.4}
and the asymptotics of the Bessel function, one obtains that
\eqref{d.9} has asymptotic behaviour $cn^{-1/2}(1+o(1))$, as $n\to\infty$.

2.
Using the above observation and an argument similar to the proof of
Proposition~\ref{thm.a.1}, one can prove that Proposition~\ref{thm.a.1}
is sharp in the following sense. Suppose that the
potential $V$ is not identically zero and is either non-negative or
non-positive. Then for some $c>0$ and all sufficiently large $n$,
$\text{(the width of the n'th eigenvalue cluster)}\geq cn^{-1/2}$.
\end{remark}

\section{Proof of the asymptotic expansions \eqref{a.7a}, \eqref{a.7}}
\label{sec.e}
We will prove the asymptotic expansion \eqref{a.7}
by expressing the eigenvalue moments $\mom_n^{(k)}$
as contour integrals of an analytic function.
Let $\G_n$ be a positively oriented circle around $\L_n$
with the radius $\B$.
First, we need estimates on the norm of the `sandwiched resolvent'
of $H$ on the contours $\G_n$.
\begin{lemma}\label{l.e.1}
For $n\to\infty$, one has
\begin{align}
\sup_{z\in \G_n}\norm{\abs{V}^{1/2}R(z)\abs{V}^{1/2}}&=O(n^{-1/2}\log n),
\label{e.2}
\\
\sup_{z\in \G_n}\norm{\abs{V}^{1/2}R(z)\abs{V}^{1/2}}_{S_2}&=O(n^{-1/4}\log n).
\label{e.3}
\end{align}
\end{lemma}
\begin{proof}
Let us prove \eqref{e.2}.
Using the estimate \eqref{d.2}, we get for $z\in\G_n$:
\begin{multline*}
\norm{\abs{V}^{1/2}R(z)\abs{V}^{1/2}}
\leq
\sum_{k=1}^\infty
\frac{\norm{\abs{V}^{1/2}P_k\abs{V}^{1/2}}}{\abs{\L_k-\l}}
\leq
\sum_{k=1}^\infty\frac{C}{\sqrt{k}\abs{\L_k-z}}
+O(n^{-1})
\\
\leq
C\int_0^{n-1}\frac{dx}{\sqrt{x}\abs{\B(2x+1)-z}}
+C\int_{n+1}^\infty\frac{dx}{\sqrt{x}\abs{\B(2x+1)-z}}
+
O(n^{-1/2})=O(n^{-1/2}\log n),
\end{multline*}
as $n\to\infty$.
The estimate \eqref{e.3} can be proven in a similar fashion by using \eqref{d.3}.
\end{proof}

The core of the proof  of
the expansions \eqref{a.7a}, \eqref{a.7} is the following Lemma.
\begin{lemma}\label{lma.e.2}
For all $k\in\Z_+$ and all $j\geq 2$, the integrals
\begin{equation}
\int_{\G_n}\Tr(VR(z))^j (z-\L_n)^{k}dz
\label{f.40}
\end{equation}
have an asymptotic expansion in integer powers of $n^{-1/2}$
as $n\to\infty$.
\end{lemma}
The proof of Lemma~\ref{lma.e.2} is given in sections~\ref{sec.f},~\ref{sec.fa}.
\begin{proof}[\textbf{Proof of the asymptotic expansions \eqref{a.7a}, \eqref{a.7}}]
First of all, note that it suffices to prove \eqref{a.7a}, \eqref{a.7}
with some \emph{complex} coefficients $\a_j^{(k)}$;
indeed, as $\mom_n^{(k)}$
are real, \emph{a posteriori} the coefficients $\a_j^{(k)}$ are easily seen to be real.
Thus, in what follows we will work with expansions with complex coefficients.

By \cite[Theorem 2.11]{AHS}, the difference of the resolvents of $H+V$ and $H$
belongs to the trace class. This enables us to define the analytic function
$$
W(z)=\Tr((H+V-z)^{-1}-(H-z)^{-1})
=
-\int_{-\infty}^\infty\frac{\xi(\l)}{(\l-z)^2}d\l,
\qquad z\in\C\setminus(\s(H)\cup\s(H+V)).
$$
The second equality in the above formula
is due to Krein's trace formula \eqref{a.3}.
Let $n$ be sufficiently large so that
$$
\supp\xi\cap\D_n\subset[\L_n-\frac{\B}{2},\L_n+\frac{\B}{2}]
$$
(see \eqref{d.11}).
Let, as above, $\G_n$ be a positively oriented circle around $\L_n$
with the radius $\B$.
Integrating $W(z)(z-\L_n)^{k+1}$ over $z$ around $\G_n$,
we obtain
\begin{multline}
-\frac{1}{2\pi i}\int_{\G_n}W(z)(z-\L_n)^{k+1}dz
=
\int_{-\infty}^\infty d\l\,\xi(\l) \,\frac{1}{2\pi i}
\int_{\G_n}\frac{(z-\L_n)^{k+1}}{(z-\l)^2}dz
\\
=
(k+1)\int_{\D_n}\xi(\l)(\l-\L_n)^{k}d\l
=
\mom_n^{(k)}.
\label{e.1}
\end{multline}
Firstly we prove the expansion \eqref{a.7} (i.e. assume $k\geq1$).
Expanding the resolvent $(H+V-z)^{-1}$ yields
\begin{equation}
W(z)=\sum_{j=1}^\infty(-1)^j\Tr[R(z)(VR(z))^j].
\label{e.4}
\end{equation}
Lemma~\ref{l.e.1} ensures that the series in \eqref{e.4}
converges absolutely for $z\in\G_n$ and large $n$.
Substituting the expansion \eqref{e.4} into  \eqref{e.1}
and subsequently integrating by parts
in each term of the series, we obtain
\begin{equation}
\mom_n^{(k)}=(k+1)
\sum_{j=k+1}^\infty\frac{(-1)^j}{j}\frac{1}{2\pi i}
\int_{\G_n}\Tr(VR(z))^j (z-\L_n)^{k}dz,
\quad k\in\N.
\label{e.6c}
\end{equation}
Here the summation starts from $j=k+1$, as for
$j\leq k$ the integrand is analytic at $z=\L_n$ and therefore
the integral vanishes.
Using  Lemma~\ref{l.e.1}, we obtain the following
estimates for the integrals in the r.h.s. of \eqref{e.6c}:
\begin{multline*}
\aabs{\int_{\G_n}\Tr(VR(z))^j (z-\L_n)^{k}dz}
\leq
\B^k\int_{\G_n}\norm{\abs{V}^{1/2}R(z)\abs{V}^{1/2}}^2_{S_2}
\norm{\abs{V}^{1/2}R(z)\abs{V}^{1/2}}^{j-2}dz
\\
\leq C\br{\frac{C\log n}{n^{1/4}}}^2\br{\frac{C\log n}{n^{1/2}}}^{j-2}.
\end{multline*}
This ensures that the series \eqref{e.6c} converges absolutely
(for sufficiently large $n$) and gives a bound for the
remainder:
\begin{equation}
\aabs{\sum_{j=N}^\infty\frac{(-1)^j}{j}\frac{1}{2\pi i}
\int_{\G_n}\Tr(VR(z))^j(z-\L_n)^kdz}
=
O\br{(\log n)^N n^{-(N-1)/2}},
\quad n\to\infty.
\label{e.7}
\end{equation}
Combining Lemma~\ref{lma.e.2} with the estimate
\eqref{e.7}, we obtain that the moments $\mom_n^{(k)}$,
$k\geq1$ have an asymptotic expansion in integer powers
of $n^{-1/2}$ as $n\to\infty$.
We also need to prove that
first several terms of the expansion vanish, so that
the expansion starts from
the term $Cn^{-k/2}$.
This can be seen as follows.
Note that for $j=k+1$ we can compute the integral in the series
\eqref{e.6c}, which gives
\begin{equation}
\mom_n^{(k)}=
\Tr(VP_n)^{k+1}+
(k+1)\sum_{j=k+2}^\infty\frac{(-1)^j}{j}\frac{1}{2\pi i}
\int_{\G_n}\Tr(VR(z))^j (z-\L_n)^{k}dz,
\quad k\in\Z_+.
\label{e.6}
\end{equation}
Lemma~\ref{l.d.1} gives
$$
\Tr(VP_n)^{k+1}\leq\norm{\abs{V}^{1/2}P_n\abs{V}^{1/2}}_{S_2}^2
\norm{\abs{V}^{1/2}P_n\abs{V}^{1/2}}^{k-1}=O(n^{-k/2}),
\quad
n\to\infty.
$$
Combining this with the estimate \eqref{e.7}
with $N=k+2$, we obtain $\mom_n^{(k)}=O(n^{-k/2})$
as $n\to\infty$.

Secondly we prove the expansion \eqref{a.7a}, i.e., the case
$k=0$. Here the only difference is that the first term in the
series \eqref{e.6c} is not well defined, as $VR(z)$ is not of the
trace class. However, this term can be written as (cf.
\eqref{e.6}) $\Tr(\abs{V}^{1/2}P_n\abs{V}^{1/2}\sign V)$, and by
the explicit form \eqref{d.1} of the integral kernel of $P_n$, the
last expression equals $\frac{\B}{2\pi}\int V(x)dx$. The rest of
the argument is the same as for $k\geq 1$.
\end{proof}

\begin{proof}[\textbf{Proof of the formula  \eqref{a.8} for $\a_3^{(0)}$}]
We will prove
\begin{equation}
\a_0^{(1)}=\frac{\B^{3/2}}{4\sqrt{2}\pi^3}
\int_{\R^2}\int_{\R^2}\frac{V(x)V(y)}{\abs{x-y}}dx\,dy.
\label{g.1}
\end{equation}
By \eqref{a.9}, this will also imply the formula \eqref{a.8}
for $\a_3^{(0)}$.

Due to \eqref{e.6} and \eqref{e.7}, we have
$$
\mom_n^{(1)}=\Tr(VP_n)^2+O((\log n)^3n^{-1}),
\quad n\to\infty,
$$
so it suffices to prove that
$$
\Tr(VP_n)^2=\frac{\a_0^{(1)}}{n^{1/2}}+o(n^{-1/2}),
\quad n\to\infty,
$$
with $\a_0^{(1)}$ given by \eqref{g.1}.
By formula \eqref{d.1} for the integral kernel of $P_n$,
we have
\begin{multline*}
\Tr(VP_n)^2=\br{\frac{B}{2\pi}}^2
\int_{\R^2\times\R^2} V(x)V(y)L_n(\tfrac{\B}{2}\abs{x-y}^2)^2
\exp(-\tfrac{\B}{2}\abs{x-y}^2)dx\, dy
\\
=
\frac{\B}{2\pi^2}\int_0^\infty  L_n(t^2)^2 e^{-t^2}h(t)tdt,
\end{multline*}
where $h\in C_0^\infty(\R)$ is given by
$$
h(t)=\int_{{\mathbb S}^1} d\o \int_{\R^2}dy
V(y)V(y+\sqrt{\tfrac{2}{\B}}t\o),
\quad t\in\R.
$$
By \eqref{d.4},
$$
\int_0^\infty  L_n(t^2)^2 e^{-t^2}h(t)tdt=
\int_0^\infty  J_0(t\sqrt{4n+2})^2h(t)tdt+O(n^{-3/4}),
\quad n\to\infty.
$$
Next, using the asymptotics of the Bessel function, we obtain
\begin{gather}
\bigl\lvert J_0(x)-\sqrt{\tfrac{2}{\pi x}}\cos(x-\tfrac{\pi}{4})
\bigr\rvert
\leq C x^{-1/2}(1+x)^{-1},\quad x>0,
\notag
\intertext{and therefore}
\bigl\lvert J_0(x)^2-\tfrac{2}{\pi x}(\cos(x-\tfrac{\pi}{4}))^2
\bigr\rvert
\leq C x^{-1}(1+x)^{-1},\quad x>0.
\label{e.8}
\end{gather}
One has
\begin{gather}
\begin{split}
\int_0^\infty\tfrac{2}{\pi t\sqrt{4n+2}}
(&\cos(t\sqrt{4n+2}-\tfrac{\pi}{4}))^2 h(t)tdt =
\tfrac{1}{\pi\sqrt{4n+2}}\int_0^\infty h(t)dt
\\
&+
\tfrac{1}{\pi\sqrt{4n+2}}\int_0^\infty
\cos(2t\sqrt{4n+2}-\tfrac{\pi}{2})h(t)dt
=
\tfrac{1}{2\pi\sqrt{n}}\int_0^\infty h(t)dt
+
o(n^{-1/2}),\quad n\to\infty;
\end{split}
\label{e.9}
\\
\int_0^\infty
\tfrac{1}{t\sqrt{4n+2}(1+t\sqrt{4n+2})}h(t)tdt
=
o(n^{-1/2}),\quad n\to\infty.
\label{e.10}
\end{gather}
Combining \eqref{e.8}, \eqref{e.9}, \eqref{e.10},
and computing the integral
$\int_0^\infty h(t)dt$ yields
formula \eqref{g.1}.
\end{proof}

\section{Analytic properties of the resolvent  $R(z)$}
\label{sec.f}
In this section, we discuss analytic properties of the integral
kernel of the resolvent $R(z)=(H-z)^{-1}$
and reduce the proof of Lemma~\ref{lma.e.2} to Lemma~\ref{l.f.1}.
Our analysis is based on the following explicit formula for this
kernel:
\begin{lemma}\label{l.kernel}
For any $z\in\C\setminus\s(H)$,
the integral kernel of the resolvent $R(z)$ of the magnetic Hamiltonian $H$
can be expressed in terms of $\G$-function
and confluent hypergeometric function $U(a,b;\zeta)$
as follows:
\begin{equation}
R(z)(x,y)=\frac{1}{4\pi}\G\br{\tfrac12-\tfrac{z}{2\B}}
U\br{\tfrac12-\tfrac{z}{2\B},1;\tfrac{\B}{2}\abs{x-y}^2}
\exp\br{-\tfrac{\B}{4}\abs{x-y}^2+i\tfrac{\B}{2}[x,y]},
\label{f.1}
\end{equation}
where $x,y\in\R^2$, $x\not=y$.
\end{lemma}
\begin{proof}
Let us employ the integral representation \cite[(13.2.5)]{AS}
for the confluent hypergeometric function
\begin{equation}
\G(a)U(a,1;\zeta)=
\int_0^\infty e^{-\zeta\tau}\tau^{a-1}(1+\tau)^{-a}d\tau,
\quad 0<\Re a<1,\quad \z>0
\label{f.2}
\end{equation}
and  the explicit formula \eqref{f.2a}
for the heat kernel \cite{AHS}
of the magnetic Hamiltonian $H$.
Substituting \eqref{f.2a} into the formula
$$
R(z)=\int_0^\infty e^{-tH}e^{tz}dt,
$$
denoting $\z=\frac{\B}{2}\abs{x-y}^2$,
making the change of variable $\tau=(\coth(\B t)-1)/2$ in the integral,
and taking into account \eqref{f.2}, one obtains \eqref{f.1}
for $-\B<\Re z<\B$. Analytic continuation in $z$
completes the argument.
\end{proof}
For the reader's convenience and ease of further reference
we start by recalling the necessary facts about the confluent hypergeometric
functions $U(a,b;\zeta)$, $M(a,b;\zeta)$.
Our main sources are \cite{AS} and \cite{Tri}.
The functions $U(a,b;\zeta)$ and $M(a,b;\zeta)$ are two linearly
independent solutions to the Kummer's equation
$$
\zeta \frac{d^2 U}{d\zeta^2}+(b-\zeta)\frac{dU}{d\zeta}-aU=0.
$$
We are only interested in the case $b=1$ (see \eqref{f.1})
or $b$ lying in a small neighbourhood
of $1$, so we assume that $\abs{b-1}<1/2$;
this will simplify our discussion.
We also assume $0<\zeta\leq R$ for some fixed $R>0$,
as we are interested in the case $\zeta=\frac{\B}{2}\abs{x-y}^2$
when both $x$ and $y$ are in $\supp V$ (see \eqref{f.40}).

The function $M(a,b;\zeta)$ is given by a convergent Taylor series
$$
M(a,b;\zeta)
=
1+\frac{a}{b}\frac{\zeta}{1!}+\frac{a(a+1)}{b(b+1)}\frac{\zeta^2}{2!}
+
\frac{a(a+1)(a+2)}{b(b+1)(b+2)}\frac{\zeta^3}{3!}+\cdots
$$
As it is readily seen from the above series, $M(a,b;\z)$ is
analytic in
$(a,b,\z)\in\C\times\{b:\abs{b-1}<1/2\}\times\C$.
For $-a\notin\Z_+$, $\z>0$, $b\not=1$,
the function $U(a,b;\zeta)$ is defined by
\begin{equation}
\G(a)U(a,b;\zeta)=\frac{\pi}{\sin (\pi b)}
\br{\frac{\G(a)}{\G(1+a-b)\G(b)}M(a,b;\zeta)
-
\z^{1-b}\frac{M(1+a-b,2-b;\zeta)}{\G(2-b)}}.
\label{f.3}
\end{equation}
We assume that $\arg\z=0$; this fixes the branch of
$\z^{1-b}$.
The function $\G(a)U(a,b;\z)$ is meromorphic in $a\in\C$ with poles at
$a=0,-1,-2,\dots$ which correspond to the Landau levels --- see
\eqref{f.1}.

The r.h.s. of \eqref{f.3} is analytic in $b$ with a removable
singularity at $b=1$; the limit as $b\to1$ is easy to compute:
\begin{equation}
\G(a)U(a,1;\z)=-2M'_b(a,1;\z)-
M'_a(a,1;\z)-(2\g+\psi(a)-\log\z)M(a,1;\z),
\label{f.4}
\end{equation}
where
$M'_a=\frac{\partial M}{\partial a}$,
$M'_b=\frac{\partial M}{\partial b}$,
 $\psi(a)=\G'(a)/\G(a)$ is the digamma function,
$\g$ is Euler's constant $\g=-\psi(1)\approx0.577$,
and $\log\z\in\R$, $\z>0$.

Using the reflection formula for the $\psi$ function,
$$
\psi(a)=\psi(1-a)-\pi\cot(\pi a),
$$
let us rearrange formula \eqref{f.4} as
\begin{align}
\G(a)U(a,1;\z)&=\wt M(a,\z)+\pi\cot(\pi a)M(a,1;\z),
\label{f.6}
\\
\wt M(a,\z)&=-2M'_b(a,1;\z)-M'_a(a,1;\z)
-(2\g+\psi(1-a)+\log\z)M(a,1;\z).
\label{f.7}
\end{align}
The function $\wt M(a,\z)$ is analytic in $a$ at the points
$-a\in\Z_+$. The singularities of $\G(a)U(a,1;\z)$ written in
the form \eqref{f.6} are easy to analyze, as they are due to
the elementary function $\cot{\pi a}$.
Incidentally, \eqref{f.6} gives a formula for the residues
of the resolvent $R(z)$; due to the identity $M(-n,1;\z)=L_n(\z)$,
this formula agrees with \eqref{d.1}.

\begin{proof}[\textbf{Proof of Lemma~\ref{lma.e.2}}]
Substituting formula \eqref{f.1} into the integrals \eqref{f.40}
and using \eqref{f.6}, we see that in order to obtain the required
asymptotic expansions, we need to analyse the asymptotics
of the integrals
\begin{equation}
\int_{\g_n}(a+n)^k(\cot\pi a)^u G(a)da,
\quad n\to\infty,
\quad k\in\Z_+,
\quad u\in\N,
\label{f.5aa}
\end{equation}
where $\g_n$ is a positively oriented circle around $a=-n$
with the radius $1/2$, and
$G(a)$ is the analytic function
\begin{equation}
G(a)=\int_{\R^{2j}}F(x_1,\dots,x_j)
\Pi_{p=1}^j \M_p(a,\tfrac{\B}{2}\abs{x_{p+1}-x_p}^2)
dx_1\cdots dx_j,
\qquad x_{j+1}\equiv x_1.
\label{f.5a}
\end{equation}
Here
$$
F(x_1\dots x_j)=V(x_1)\cdots V(x_j)
\exp\br{-\tfrac{\B}{4}\sum_{p=1}^j\abs{x_{p+1}-x_p}^2
+i\tfrac{\B}{2}\sum_{p=1}^j[x_p,x_{p+1}]},
\quad x_{j+1}=x_1,
$$
each of the functions $\M_p(a,\z)$
is either $M(a,1;\z)$ or $\wt M(a,\z)$
and at least one of the functions $\M_p(a,\z)$
is $M(a,1;\z)$.

Applying the residue formula to the integral \eqref{f.5aa},
we see that the required statement follows from
 Lemma~\ref{l.f.0} below.
\end{proof}

\begin{lemma}\label{l.f.0}
For any $F\in C_0^\infty(\R^{2j})$, let $G(a)$ be given by
\eqref{f.5a}, where each of the functions $\M_p(a,\z)$ is either
$M(a,1;\z)$ or $\wt M(a,\z)$ and at least one of the functions
$\M_p(a,\z)$ is $M(a,1;\z)$. Then the function $G(a)$ and all of
its derivatives $G^{(s)}(a)$, $s\geq1$, admit asymptotic
expansions in integer powers of $\abs{a}^{-1/2}$ as $a\to-\infty$,
$a\in\R$.
\end{lemma}
\begin{proof}
The proof is based on  using suitable
asymptotic expansions of the functions
$M(a,1;\z)$ and $\wt M(a,\z)$
in terms of Bessel functions and on application
of Lemma~\ref{l.f.1}.

(i) First consider the special case
when  $\M_p(a,\z)=M(a,1;\z)$ for all $p$ in \eqref{f.5a}.
Our main tool will be a convergent expansion of
$M(a,b;\z)$ in terms of Bessel functions
due to Tricomi \cite{Tri}.
For our purposes it suffices to consider the following range of
parameters:
\begin{equation}
\Re a\leq -1,\quad \abs{\Im a}\leq 1,
\quad \abs{b-1}\leq 1/2,
\quad 0<\z\leq R
\label{f.10}
\end{equation}
for some fixed $R>0$.
The expansion of \cite{Tri} (see also \cite[(13.3.7)]{AS}) reads:
\begin{equation}
M(a,b;\z)=\G(b)e^{\z/2}\br{\frac{(b-2a)\z}{2}}^{(1-b)/2}
\sum_{m=0}^\infty \br{\frac{\z}{2b-4a}}^{m/2}A_m
J_{b-1+m}(\sqrt{(2b-4a)\z}),
\label{f.8}
\end{equation}
where $J_{b-1+m}$ are Bessel functions and $A_m=A_m(a,b)$
are the coefficients
in the Taylor expansion
\begin{equation}
f(z)=\br{\frac{e^z}{1+z}}^b
\br{e^{2z}\frac{1-z}{1+z}}^{-a}
=
\sum_{m=0}^\infty A_m z^m,
\quad \abs{z}<1.
\label{f.9}
\end{equation}
Note that $\Re (b-2a)\geq 1$ and the principal values of
$\br{\frac{(b-2a)\z}{2}}^{(1-b)/2}$ and
$\sqrt{(2b-4a)}$ are taken in \eqref{f.8}.
Due to a fast decay of the Bessel function $J_\nu(z)$ for $\nu\to\infty$,
the series \eqref{f.8} converges absolutely for
the range of parameters \eqref{f.10}.
Take $b=1$ and
for a given $N\in\N$, write Tricomi's expansion \eqref{f.8} as
\begin{equation}
M(a,1;\z)=M_N^{(0)}(a,\z)+M^{(1)}_N(a,\z),
\label{f.51}
\end{equation}
where
$$
M_N^{(0)}(a,\z)=e^{\z/2}
\sum_{m=0}^{N-1} A_m(a,1)\br{\frac{\z}{2-4a}}^{m/2}
J_{m}(\sqrt{(2-4a)\z}).
$$
Let us recall the argument of \cite{Tri} which gives
the estimate for $M_N^{(1)}(a,\z)$.
By inspecting the integral representation for the Bessel
function, one obtains a uniform estimate
\begin{equation}
\abs{J_{m}(\sqrt{(2-4a)\z})}\leq C \quad \text{ for }\quad
m\in\Z_+,\quad\Re a\leq -1,\quad \abs{\Im a}\leq 1,\quad
\quad 0< \z\leq R.
\label{f.10a}
\end{equation}
Next,  one needs to estimate for the coefficients $A_m$
of the expansion \eqref{f.8}.
Applying Cauchy's theorem to the Taylor expansion \eqref{f.9}
yields
$$
\abs{A_m}\leq r^{-m}\max_{\abs{z}=r}\abs{f(z)}
$$
for any $r\in(0,1)$. Note that
$$
f(z)=\br{\frac{e^z}{1+z}}^b(1+O(z^3))^{-a},
\quad z\to0.
$$
Choosing $r=\abs{a}^{-\frac{5}{12}}$
(one can take $r=\abs{a}^{-\delta}$ for any $1/3<\delta<1/2$),
one obtains for all sufficiently large $\abs{a}$:
\begin{equation}
\abs{A_m}\leq C\abs{a}^{\frac{5}{12} m}
(1+C\abs{a}^{-\frac{15}{12}})^{-\text{Re}\, a}
\leq
2C\abs{a}^{\frac{5}{12}m}.
\label{f.11}
\end{equation}
Combining \eqref{f.10a} and \eqref{f.11} gives
for all sufficiently large $\abs{a}$:
\begin{equation}
\abs{M_N^{(1)}(a,\z)}\leq
C\sum_{m=N}^\infty\abs{a}^{\frac{5}{12} m}
\aabs{\frac{R}{2-4a}}^{m/2}
=O(\abs{a}^{-\frac1{12}N}),
\quad\Re a\to-\infty,
\quad \abs{\Im a}\leq1.
\label{f.12}
\end{equation}
In the same way, the estimates \eqref{f.10a} and \eqref{f.11}
show that
\begin{equation}
\abs{M_N^{(0)}(a,\z)}=O(1),\qquad \Re a\to-\infty,
\quad \abs{\Im a}\leq1.
\label{f.53}
\end{equation}
Substituting \eqref{f.51} into \eqref{f.5a}, we obtain
$$
G(a)=G_N^{(0)}(a)+G_N^{(1)}(a),
$$
where
$$
G_N^{(0)}(a)=\int_{\R^{2j}}F(x_1,\dots,x_j)
\Pi_{p=1}^j M_N^{(0)}(a,\tfrac{\B}{2}\abs{x_{p+1}-x_p}^2)
dx_1\cdots dx_j,
\qquad x_{j+1}\equiv x_1.
$$
By \eqref{f.53}, \eqref{f.12}, we get
$$
G_N^{(1)}(a)=O(\abs{a}^{-\frac{1}{12}N}),
\qquad \Re a\to-\infty,
\quad \abs{\Im a}\leq 1.
$$
By Cauchy's formula for the derivatives, this entails
$$
\br{\frac{d}{da}}^s G_N^{(1)}(a)=
O(\abs{a}^{-\frac{1}{12}N}),
\quad a\to-\infty,\quad a\in\R.
$$
As $N$ can be taken arbitrary large, we see that
it suffices to prove that for any $N>0$, all derivatives
$\br{\frac{d}{da}}^s G_N^{(0)}(a)$, $s\in\Z_+$,
have an asymptotic expansion for $a\to-\infty$, $a\in\R$.

{}From \eqref{f.9} it follows that
the coefficients $A_m(a,b)$ are polynomials in $a$
and $b$.
This observation reduces the problem to justifying the
asymptotic expansion of the integral \eqref{f.5a},
where each of the functions $\M_p(a,\z)$ is
$\z^{m/2} J_m(\sqrt{(2-4a)\z})$ with some $m\in\Z_+$.
Such an expansion is provided by Lemma~\ref{l.f.1}.

(ii) Consider the general case.
First let us obtain an expansion for $\wt M(a,\z)$
similar to \eqref{f.8}.
Substituting \eqref{f.8} into the r.h.s. of \eqref{f.7},
after a rearrangement we obtain
\begin{multline}
\wt M(a,\z)
=
-(\psi(1-a)-\log (\frac12-a))
e^{\z/2}\sum_{m=0}^\infty A_m\br{\frac{\z}{2-4a}}^{m/2}
J_{m}(\sqrt{(2-4a)\z})
\\
-
2 e^{\z/2}\sum_{m=0}^\infty A_m
\br{\frac{\z}{2-4a}}^{m/2}
\dot J_{m}(\sqrt{(2-4a)\z})
-
e^{\z/2}\sum_{m=0}^\infty B_m
\br{\frac{\z}{2-4a}}^{m/2}
J_{m}(\sqrt{(2-4a)\z}),
\label{f.15}
\end{multline}
where
$$
B_m=\br{2\frac{\partial A_m}{\partial b}+\frac{\partial A_m}{\partial a}}
\mid_{b=1},
\quad
\dot J_m(z)=\frac{\partial J_\nu(z)}{\partial\nu}\mid_{\nu=m}.
$$
A fast decay of $J_m$ and $\dot J_m$ as $m\to\infty$
ensures convergence of the
series and validates differentiation with respect to  $a$ and $b$.

Using the expansion \eqref{f.15}, we can complete the
argument by following the same steps as in part (i) of the proof.
First we  need to obtain estimates
for the remainders of the series in the r.h.s. of
\eqref{f.15} in the strip $\Re a\leq -1$, $\abs{\Im a}\leq 1$
(cf. \eqref{f.12}).
The estimate for the remainder term of the first series in the r.h.s.
of \eqref{f.15} is provided by \eqref{f.12}.
The estimate for the second series in the r.h.s.
of \eqref{f.15}
is obtained in exactly the same way by using the estimates
$\abs{\dot J_0(\sqrt{(2-4a)\z})}\leq
C+C\abs{\log\abs{(2-4a)\z}}$ and (see \cite[(9.1.22)]{AS})
$$
\quad
\abs{\dot J_m(\sqrt{(2-4a)\z})}\leq C,
\quad m\in\N,\quad \Re a\leq -1,
\quad \abs{\Im a}\leq 1,
\quad 0<\z\leq R
$$
instead of \eqref{f.10a}.
In order to estimate the remainder term of the third series,
we need an estimate on the coefficients $B_m$.
The coefficients $B_m$ are readily seen to be Taylor coefficients
of the function (cf. \eqref{f.9})
$$
g(z)=\br{2\frac{\partial f}{\partial b}(z)+\frac{\partial f}{\partial a}(z)}
=f(z)\log (1-z^2)^{-1/2}
$$
which similarly to \eqref{f.11} gives
$$
\abs{B_m}\leq C\abs{a}^{\frac{5}{12} m}.
$$
This gives the analogue of the estimate \eqref{f.12} for the third series in
the r.h.s. of \eqref{f.15}.

Next, the function
$\psi(1-a)-\log (\frac12-a)$ in \eqref{f.15}
admits asymptotic expansion
in integer powers of $a^{-1}$ as $\Re a\to-\infty$
(see \cite[6.3.18]{AS}).
Thus, we have reduced the problem to justifying an asymptotic
expansion of the integral \eqref{f.5a},
where each of the functions $\M_p(a,\z)$ is either
$\z^{m/2} J_m(\sqrt{(2-4a)\z})$
or $\z^{m/2} \dot J_m(\sqrt{(2-4a)\z})$
and at least one of the functions $\M_p(a,\z)$ is
$\z^{m/2} J_m(\sqrt{(2-4a)\z})$.
Finally, the formula \cite[(9.1.66)]{AS}
$$
(\z/2)^m\dot J_m(\z)=\frac{\pi}2(\z/2)^m Y_m(\z)
+\frac{m!}{2}\sum_{k=0}^{m-1}
\frac{(\z/2)^k}{(m-k)!k!} J_k(\z)
$$
reduces the problem to Lemma~\ref{l.f.1}.
\end{proof}

\section{Asymptotics of integrals containing Bessel functions}
\label{sec.fa}

\begin{lemma}\label{l.f.1}
Let $F\in C_0^\infty(\R^{2j})$. Define a function $G_1(\aa)$,
$\aa>0$, by
\begin{equation}
G_1(\aa):=
\int_{\R^{2j}}F(x_1,\dots x_j)
\prod_{p=1}^j
\J_{m_p}(\aa\abs{x_p-x_{p+1}})\abs{x_p-x_{p+1}}^{m_p}
dx_1\cdots dx_j,
\quad x_{j+1}=x_j
\label{f.16}
\end{equation}
where each of the functions
$\J_{m_p}(\z)$ is either $J_{m_p}(\z)$ or
$Y_{m_p}(\z)$ with some $m_p\in\Z_+$, and at least one
of the functions $\J_{m_p}(\z)$ is $J_{m_p}(\z)$.
Then the function $G_1(\aa)$ and all of its derivatives
$G_1^{(s)}(\aa)$, $s\geq1$, have a complete asymptotic
expansion in integer powers of $\aa^{-1}$ for $\aa\to+\infty$.
\end{lemma}
\begin{proof}
Recall the identity \cite[(9.1.27)]{AS}
$$
\z^{\nu+1} \J_{\nu+1}(\z)=2\nu\z^\nu J_\nu(\z)
-\z^2(\z^{\nu-1}J_{\nu-1}(\z)),
\quad \J_\nu=J_\nu\text{ or } \J_{\nu}=Y_\nu.
$$
This identity allows us to reduce the problem to the case
when all the indices $m_p$ in the integral \eqref{f.16}
are $0$ or $1$.
Next, assume for the convenience of notation
that the first $l$ functions $\J$ in the integral \eqref{f.16}
are the Neumann functions $Y$,
 and the remaining $j-l$
functions are the Bessel functions $J$.
For $\aa_1>0$, $\aa_2>0$, \dots, $\aa_j>0$
define
\begin{multline}
G_2(\aa_1,\dots,\aa_j)
\\
=\int_{\R^{2j}}F(x_1,\dots x_j)
\prod_{p=1}^l
Y_{0}(\aa_p\abs{x_p-x_{p+1}})
\prod_{p=l+1}^j
J_{0}(\aa_p\abs{x_p-x_{p+1}})
dx_1\cdots dx_j,
\quad
x_{j+1}\equiv x_1.
\label{f.17}
\end{multline}
Formulae
$$
\frac{d J_0(\aa\abs{x})}{d\aa}=
-\abs{x} J_1(\aa\abs{x}),
\quad
\frac{d Y_0(\aa\abs{x})}{d\aa}=
-\abs{x} Y_1(\aa\abs{x})
$$
show that it suffices to obtain an asymptotic expansion of the
functions
\begin{equation}
\br{\tfrac{\partial}{\partial \aa_1}}^{\b_1}\cdots
\br{\tfrac{\partial}{\partial \aa_j}}^{\b_j}
 G_2(\aa_1,\dots \aa_j)\mid_{\aa_1=\cdots=\aa_j=\aa},
\qquad \b_p\in\Z_+
\label{f.18}
\end{equation}
for $\aa\to\infty$.
We shall obtain an asymptotic expansion for the function
$G_3(\aa):=G_2(\aa,\dots,\aa)$ as $\aa\to\infty$.
{}From the construction it will be clear that the derivatives
\eqref{f.18} can be dealt with in  the same way.
Let us make a change of variables in the integral \eqref{f.17}.
Denote
\begin{gather*}
y_p=x_p-x_{p+1},\quad p=1,\dots, j-1,
\quad z=x_1+\cdots+x_j,
\\
F_1(y)=\int_{\R^2} F(x_1\dots x_j)dz,\qquad F_1\in C_0^\infty(\R^{2j-2}).
\end{gather*}
Then we obtain
\begin{equation}
G_3(\aa)=\int_{\R^{2j-2}}
F_1(y) J_0(\aa\abs{y_1+\dots+y_{j-1}})
\prod_{p=1}^l Y_0(\aa\abs{y_p})
\prod_{p=l+1}^{j-1} J_0(\aa\abs{y_p})
dy.
\label{f.19}
\end{equation}
Next, we use the formulae
$$
J_0(\aa\abs{y})=\frac{1}{\pi}\int_{\R^2} e^{i\aa uy}\d(u^2-1)du,
\qquad
Y_0(\aa\abs{y})=-\frac{1}{\pi^2}\vp \int_{\R^2} \frac{e^{i\aa uy}}{u^2-1}du.
$$
Substituting these formulae into \eqref{f.19}, we obtain
$$
G_3(\aa)
=
\vp\int_{\R^{2j}}
\frac{\d(u_{l+1}^2-1)\cdots\d(u_{j}^2-1)}{(u_1^2-1)\cdots(u_l^2-1)}
F_2(\aa(u_1-u_j,u_2-u_j,\dots,u_{j-1}-u_j))
du_1\dots du_j,
$$
where $F_2$ is (up to a multiplicative constant) the Fourier transform of $F_1$.
Note that $F_2$ belongs to the Schwartz class ${\cal S}(\R^{2j-2})$.

In order to simplify the last integral, we introduce some notation.
Let us use the polar coordinates $u_i=r_i^{1/2}\oo_i$,
where $\oo_i=(\cos\o_i,\sin\o_i)\in\R^2$, $\o_i\in\T=\R/2\pi\Z$.
Denote also $\o=(\o_1,\dots,\o_{j-1})\in\T^{j-1}$,
$r=(r_1,\dots,r_l)\in\R^l_+$.
Define the function
$$
f(r,\o,\o_j)=
(r_1^{1/2}\oo_1-\oo_j,\dots,r_l^{1/2}\oo_l-\oo_j,
\oo_{l+1}-\oo_j,\dots,\oo_{j-1}-\oo_j)\in\R^{2j-2}.
$$
With this notation, we have
\begin{equation}
G_3(\aa)=2^{-j}\int_\T d\o_j
\int_{\T^{j-1}}d\o \text{ v.p.}\int_{\R_+^l}dr
\frac{F_2(\aa f(r,\o,\o_j))}{(r_1-1)\cdots(r_l-1)}.
\label{f.20}
\end{equation}
Note that
$$
f(r,\o,\o_j)=0\Leftrightarrow
\left( r=(1,\dots,1)\in\R_+^l \text{ and }
\o=(\o_j,\dots,\o_j)\in\T^{j-1} \right)
$$
and $\rank f'(r,\o,\o_j)=l+j-1$ at the point $r=(1,\dots,1)$, $\o=(\o_j,\dots,\o_j)$.
Let us show that only an arbitrary small neighbourhood of the point
$r=(1,\dots,1)$, $\o=(\o_j,\dots,\o_j)$ gives contribution to the
asymptotics of the integral \eqref{f.20}.
First  recall some estimates for the principal value integrals.
Let $\delta>0$ and $\phi\in C^\infty(-\d,\d)$. Then
$$
\vp \int_{-\d}^\d\frac{\phi(x)}{x}dx=
\int_{-\d}^\d\phi'(x)\log\abs{\d/x}dx,
$$
and using the Cauchy-Schwartz inequality, one obtains
\begin{equation}
\aabs{\vp \int_{-\d}^\d\frac{\phi(x)}{x}dx}\leq
C(\d)\norm{\phi'}_{L^2(-\d,\d)}.
\label{f.21}
\end{equation}
Similarly, for $\phi\in C^\infty([-\d,\d]^l)$,
\begin{equation}
\aabs{\vp \int_{(-\d,\d)^l} \frac{\phi(x)}{x_1\cdots x_l}dx}\leq
C(\d)
\nnorm{\frac{\partial^l\phi}{\partial x_1\cdots \partial x_l}}_{L^2((-\d,\d)^l)}.
\label{f.22}
\end{equation}
Denote
$U=(1-\e,1+\e)^l\times (\o_j-\e,\o_j+\e)^{j-1}\subset \R_+^l\times\T^{j-1}$
where $\e>0$ is sufficiently small.
Let us  show that
\begin{equation}
G_3(\aa)=2^{-j}\int_\T d\o_j \vp\int_U dr d\o
\frac{F_2(\aa f(r,\o,\o_j))}{(r_1-1)\cdots (r_l-1)}
+O(\aa^{-\infty}).
\label{f.22a}
\end{equation}
For simplicity consider the case $j=2$, $l=1$.
Then for $\aa\to\infty$ one has
$$
\aabs{\int_{\T}d\o_2\int_{\T} d\o_1 \int_{\R_+\setminus(1-\e,1+\e)}
\frac{dr_1}{r_1-1} F_2(\aa f(r_1,\o_1,\o_2))}
\leq
C\sup\limits_{\abs{r_1-1}>\e}\abs{F_2(\aa f(r_1,\o_1,\o_2))}
=O(\aa^{-\infty}),
$$
and also by \eqref{f.21}
\begin{multline*}
\aabs{\int_{\T}d\o_2\, \vp\int_{1-\e}^{1+\e}
\frac{dr_1}{r_1-1}\int_{\T\setminus(\o_2-\e,\o_2+\e)}
d\o_1 F_2(\aa f(r_1,\o_1,\o_2))}^2
\\
\leq
C\int_{\T}d\o_2\int_{\T\setminus(\o_2-\e,\o_2+\e)}
d\o_1
 \int_{1-\e}^{1+\e}dr_1
\aabs{\frac{\partial F_2(\aa f(r_1,\o_1,\o_2))}{\partial r_1}}^2
=O(\aa^{-\infty}),
\end{multline*}
as $F_2$ is the Schwartz class function.

Thus, it suffices to prove an asymptotic expansion of the integral
in the r.h.s. of \eqref{f.22a}.
The asymptotic expansion of the integral over $(r,\o)$
is provided by Lemma~\ref{l.f.2} below.
It remains to prove that the asymptotic expansion given by
Lemma~\ref{l.f.2} is uniform in $\o_j$. In order to show this,
let us introduce the matrix
$$
\psi_\tau=
\begin{pmatrix}\cos\tau &-\sin\tau
\\
\sin\tau&\cos\tau
\end{pmatrix},\quad \tau\in\T,
$$
and the vector $e=(1,0)\in\R^2$.
We rewrite the function $f(r,\o,\o_j)$ in the form
\begin{multline}
f(r,\o,\o_j)=\Psi_{\o_j}\br{(r_1^{1/2}\psi_{\nu_1}-1)e,\dots,
(r_l^{1/2}\psi_{\nu_l}-1)e,(\psi_{\nu_{l+1}}-1)e,\dots,
(\psi_{\nu_{j-1}}-1)e}
\\
=\Psi_{\o_j}f(r,\nu,0),
\quad \nu_p=\o_p-\o_j,
\quad \nu=(\nu_1,\dots,\nu_{j-1})\in\T^{j-1},
\label{f.28}
\end{multline}
where $\Psi_\tau=\diag(\psi_\tau,\dots, \psi_\tau)$ is a matrix in $\R^{2j-2}$.
Substituting \eqref{f.28} into \eqref{f.22a}, we obtain
$$
G_3(\aa)=2^{-j}\int_\T d\o_j \vp\int_{U_1} dr d\nu
\frac{F_2(\aa \Psi_{\o_j}f(r,\nu,0))}{(r_1-1)\cdots (r_l-1)},
\quad
U_1=(1-\e,1+\e)^l\times (-\e,\e)^{j-1}\subset \R_+^l\times\T^{j-1}.
$$
{}From the last formula and the proof of Lemma~\ref{l.f.2}
it is clear that the expansion given by
Lemma~\ref{l.f.2} is uniform in $\o_j$;
integrating this expansion over $\o_j$, we get
the required expansion for $G_3(\aa)$.
\end{proof}
\begin{lemma}\label{l.f.2}
Let $f\in C^\infty(\R^m,\R^n)$, $m\leq n$, and suppose that
$f(0)=0$, $\rank f'(0)=m$.
Then for any sufficiently small open neighbourhood of zero
$U\subset \R^m$, any $F\in {\cal S}(\R^n)$, and $l\in\{0,1,\dots,m\}$,
the integral
\begin{equation}
I(\aa)=\vp\int_U \frac{F(\aa f(x))}{x_1\cdots x_l}dx
\label{I}
\end{equation}
has a complete asymptotic expansion
$$
I(\aa)=\aa^{l-m}(c_0+c_1\aa^{-1}+c_2\aa^{-2}+\cdots),
\quad \aa\to\infty.
$$
\end{lemma}
\begin{proof}
Choose $U$ sufficiently small so that
\begin{equation}
\abs{f(x)-f'(0)x}\leq\frac12\abs{f'(0)x},
\quad \forall x\in U.
\label{f.23}
\end{equation}
For a given $N\in\N$, $N>l$, let us prove that
\begin{equation}
I(\aa)=\aa^{l-m}\sum_{i=0}^{N-m-1}c_i \aa^{-i}+O(\aa^{l-N}),
\quad \aa\to\infty.
\label{f.23a}
\end{equation}
By Taylor's formula for $f(x)$ and $F(\aa f(x))$,
\begin{gather}
f(x)=f'(0)x+f_2(x), \qquad
f_2(x)=\sum_{s=1}^N
\frac1{s!}f^{(s)}(0)x^s+f_N(x),
\label{f.23b}
\\
F(\aa f(x))=F(\aa f'(0)x+\aa f_2(x))=
\sum_{q=0}^N F^{(q)}(\aa f'(0)x)(\aa f_2(x))^q
+F_N(x,\aa),
\label{f.24}
\\
F_N(x,\aa)=\frac{1}{N!}\int_0^1
(1-\tau)^N\br{\tfrac{d}{d\tau}}^{N+1}F(\aa f'(0)x+\aa \tau f_2(x))d\tau.
\label{f.25}
\end{gather}
Here we use simplified notation; $f^{(s)}(0)x^s$ stands for the
polylinear form of the $s$'th differential of $f$ at zero, etc.

Substituting \eqref{f.23b} into \eqref{f.24}
and collecting the terms that contain
and that do not contain $f_N(x)$ into two
different sums, we can write
\begin{equation}
F(\aa f(x))=\sum_{0\leq q\leq N}\wt F_q(\aa x) P_q(x)+
\sum_{q\geq N+1}\wt F_q(\aa x) g_q(x)+
F_N(x,\aa).
\label{f.26}
\end{equation}
Here both sums over $q$ are finite, $\wt F_q\in{\cal S}(\R^m)$
are obtained from various components of derivatives of $F$,
$P_q(x)$ are polynomials in $x$ of degree $q$, and
$g_q\in C^\infty(U)$ are functions, satisfying the polynomial type
estimates
\begin{equation}
\aabs{
\br{\tfrac{\partial}{\partial x_1}}^{\b_1}\cdots
\br{\tfrac{\partial}{\partial x_l}}^{\b_l}
g_q(x)
}\leq
C_\b\abs{x}^{q-\abs{\b}},
\qquad \abs{\b}=\b_1+\dots+\b_l\leq q,
\quad x\in U.
\label{f.27}
\end{equation}
Consider the terms obtained by substitution of the r.h.s. of
\eqref{f.26} into the integral \eqref{I}.
First, using the estimate \eqref{f.22} and the fact that
$\wt F_q$ is a Schwartz class function, we obtain
\begin{multline*}
\vp\int_U\frac{\wt F_q(\aa x)}{x_1\cdots x_l}P_q(x)dx=
\vp\int_{\R^m}\frac{\wt F_q(\aa x)}{x_1\cdots x_l}P_q(x)dx
+O(\aa^{-\infty})
\\
=\aa^{l-q-m}
\vp\int_{\R^m}\frac{\wt F_q(x)}{x_1\cdots x_l}P_q(x)dx
+O(\aa^{-\infty}),
\qquad \aa\to\infty.
\end{multline*}
So, these terms will give contribution to the asymptotics
\eqref{f.23a}.

Next, consider the terms obtained by substitution of the second sum
in \eqref{f.26} into the integral \eqref{I}.
Using \eqref{f.27}, we obtain the estimate
$$
\nnorm{\frac{\partial^l(\wt F_q(\aa x)g_q(x))}
{\partial x_1\cdots\partial x_l}}_{L^2(U)}\leq
C\aa^{l-\frac{m}{2}-q}.
$$
By \eqref{f.22}, it follows that all the corresponding integrals are
$O(\aa^{l-N})$ as $\aa\to\infty$.

Finally, consider the term $F_N(x,\aa)$.
By \eqref{f.23}, we obtain for some $c>0$:
$$
\abs{f'(0)x+\tau f_2(x)}\geq \frac12\abs{f'(0)x}\geq c\abs{x},
\quad x\in U, \quad \tau\in(0,1).
$$
Using this fact, we obtain
$$
\abs{F_N(x,\aa)}\leq
C\aa^{N+1}\sup_{\abs{y}\geq\aa c\abs{x}}
\abs{F^{(N+1)}(y)f_2(x)^{N+1}}
$$
and therefore
$$
\norm{F_N(\cdot,\aa)}_{L^2(U)}\leq
O(\aa^{-\frac{m}{2}-N-1}),
\quad \aa\to+\infty.
$$
Similarly, one can prove the estimate
$$
\nnorm{\frac{\partial^l F_N(x,\aa)}{\partial x_1\cdots\partial x_l}}_{L^2(U)}
\leq
\aa^{l-\frac{m}{2}-N-1}.
$$
By \eqref{f.22}, it follows that the integral of $F_N$
is $O(\aa^{l-N-1})$ and
will only give contribution to the remainder term in
\eqref{f.23a}.
\end{proof}

\end{document}